\documentclass[10pt,a4paper]{amsart}

\usepackage[utf8]{inputenc}
\usepackage{t1enc}

\usepackage{amssymb,amsfonts}
\usepackage{mathtools}

\usepackage{xcolor}

\usepackage[margin=2.7cm,marginpar=2cm]{geometry}

\usepackage[unicode]{hyperref}

\usepackage[capitalize,nameinlink]{cleveref}
\numberwithin{equation}{section}

\def\equationautorefname~#1\null{Equation~(#1)\null}

\def\={\,=\,}

\usepackage{thmtools}

\makeatletter
\def\@endtheorem{\endtrivlist}
\makeatother

\declaretheorem[
style=plain,
name=Theorem,
numbered=yes,
refname={Theorem,Theorems},
Refname={Theorem,Theorems}
]{Thm}

\declaretheorem[
style=plain,
name=Conjecture,
numberlike=Thm,
refname={Conjecture,Conjectures},
Refname={Conjecture,Conjectures}
]{Conj}

\declaretheorem[
style=definition,
name=Example,
numberlike=Thm,
refname={Example,Examples},
Refname={Example,Examples}
]{Eg}

\declaretheorem[
style=plain,
name=Lemma,
numberlike=Thm,
refname={Lemma,Lemmas},
Refname={Lemma,Lemmas }
]{Lem}

\newcommand{\dd}{\mathrm{d}}

\makeatletter
\@namedef{subjclassname@2020}{%
	\textup{2020} Mathematics Subject Classification}
\makeatother

\usepackage[normalem]{ulem}

\usepackage{shuffle}

\usepackage{bbm}

\DeclareMathOperator{\id}{id}

\newcommand{\mot}{\mathfrak{m}}
\newcommand{\lmot}{\mathfrak{L}}

\let\epsilon\varepsilon

\usepackage{stackengine}
 \newcommand\brabar{\smash{\raisebox{-0.5\height-0.5\depth}{\scalebox{0.6}[0.3]{$($}}\raisebox{-\height}{$\overline{\phantom{1}}$}\raisebox{-0.5\height-0.5\depth}{\scalebox{0.6}[0.3]{$)$}}}}

\begin{document}
	
	\title[Xu's Conjecture and depths of Galois descents]{On a proof of Xu's Conjecture \\ and depths of Galois descents}
	\date{August 4, 2025}
	
       \author[Steven Charlton]{Steven Charlton}
		\address{Max Planck Institute for Mathematics \\ Vivatsgasse 7 \\ 53111 Bonn \\ Germany}
	\email{charlton@mpim-bonn.mpg.de}

	\subjclass[2020]{11M32}
	
	\keywords{Multiple zeta values, Euler sums, motivic Galois descent, (alternating) double zeta values}

	\begin{abstract} 
		We give an explicit formula proving Xu's Conjecture on alternating double zeta values.  We also discuss the limitations of Glanois' motivic Galois descent criterion in this case, as it cannot specify the depth of the descent.
	\end{abstract}

	\leavevmode
	\maketitle
	
	\vspace{-1em}	
	
	\section{Introduction}
	
	For integers \( k_1,\ldots,k_d \geq 1 \), and signs \( \epsilon_1,\ldots,\epsilon_d \in \{ \pm 1 \} \), the alternating multiple zeta values (MZV's) are defined by
	\[
		\zeta(k_1,\ldots,k_d; \epsilon_1,\ldots,\epsilon_d) = \sum_{0 < n_1 < \cdots < n_d} \frac{\epsilon_1^{n_1} \cdots \epsilon_d^{n_d}}{n_1^{k_1} \cdots n_d^{k_d}} \,
	\]
	For convergence, we required \( (k_d,\epsilon_d) \neq (1,1) \).  It is convenient to decorate \( k_j \) with a bar when the corresponding \( \epsilon_j = -1 \), so that for example \( \zeta(1, \overline{2}) \coloneqq \zeta(1, 2; 1, -1) \).   When all \( \epsilon_i = 1 \), we recover the non-alternating (or classical) MZV's 
	\[ 
	\zeta(k_1,\ldots,k_d) = \sum_{0 < n_1 < \cdots < n_d} \frac{1}{n_1^{k_1} \cdots n_d^{k_d}} \,.
	\]
	For a given alternating MZV, we define the \emph{depth} to be \( d \) (the number of indices) and the \emph{weight} to be \( k_1 + \cdots + k_d \) (the sum of the indices).  We will also refer to the case \( d = 2 \) as \emph{double} zeta values. 
	
	Multiple zeta values and their alternating versions have appeared in many contexts, including a prominent connection with knot theory and high-energy physics \cite{knotsKreimer}.  The question of whether certain linear combinations of alternating MZV's descend to evaluations in terms of non-alternating MZV's is of much interest.  Broadhurst \cite{broadhurstEnum} already identified several cases where alternating MZV's could be reduced to linear combinations of non-alternating MZV's, eventually terming such examples \emph{honorary MZV's}.  In studying explicit families of Euler sums \cite{XuWt10} in weight \( \leq 10 \) (another name for alternating MZV's), Xu conjectured the following.
	
	\begin{Conj}[Xu, Conjecture 7.2, \cite{XuWt10}] \label{conj:xu}
		For positive integers \( k \), the combination
		\[
		\zeta(2, \overline{2k}) + 2k \zeta(1, \overline{2k+1}) 
		\]
		can be expressed in terms of non-alternating double zeta values.
	\end{Conj}

	 Glanois \cite[Corollary 5.1.3]{GlanoisThesis16}, \cite[Theorem 3.8]{GlanoisBasis16} gave a \emph{motivic} criterion to determine when such a descent exists.  Recall, the motivic alternating MZV's \( \zeta^\mot(\overset{\brabar}{k_1}, \ldots, \overset{\brabar}{k_d}) \) are algebraically defined versions of the alternating MZV's, which live in a (weight-)graded connected Hopf algebra comodule \( \mathcal{H}^2 \), with action \( \Delta \colon \mathcal{H}^2 \to \mathcal{A}^2 \otimes \mathcal{H}^2\), where \( \mathcal{A}^2 = \mathcal{H}^2 / \zeta^\mot(2) \mathcal{H}^2 \).  Let \( \mathcal{L}^2 = \mathcal{A}^2 / \mathcal{A}^2_{>0} \cdot \mathcal{A}^2_{>0} \).  The motivic derivations \( D_r \colon \mathcal{H}^2 \to \mathcal{L}^2_r \otimes \mathcal{H}^2 \) are defined by \( D_r \coloneqq (\pi_r \otimes \id) \circ (\Delta - 1 \otimes \id)  \), where \( \pi_r \) is the projection \( \mathcal{A}^2 \to \mathcal{L} \to \mathcal{L}_r \), to the weight \( r \) component of \( \mathcal{L} \).  The motivic iterated integrals \( I^\mot(a_0; a_1,\ldots,a_w; a_{w+1}) \), lift the classical iterated integrals, see \eqref{eqn:int}.  We then have
	\[
		\zeta^\mot(k_1,\ldots,k_d;\epsilon_1,\ldots,\epsilon_d) = (-1)^d I^\mot(0; \smash{\overbrace{\eta_1, 0, \ldots, 0}^{k_1}}, \ldots, \smash{\overbrace{\eta_d, 0, \ldots, 0}^{k_d}};  1) \,, \qquad \eta_j = \textstyle\prod_{i=j}^d \epsilon_i \,,
	\]
	and an explicit combinational formula for \( D_r \) is given by
	\begin{align*}
		&D_r I^\mot(a_0; a_1,\ldots,a_w; a_{w+1}) \\
		&= \sum_{p=0}^{w-r} I^\lmot(a_p; a_{p+1}, \ldots, a_{p+r}; a_{p+r+1}) \otimes I^\mot(a_0; a_1,\ldots,a_p , a_{p+r+1}, \ldots, a_w; a_{w+1}) \,.
	\end{align*}
	The motivic non-alternating MZV's form a subspace \( \mathcal{H}^1 \) of \( \mathcal{H}^2 \), to which the coaction and derivations descent; set \( \mathcal{A}^1 = \mathcal{H}^1 / \zeta^\mot(2) \mathcal{H}^1 \), and \( \mathcal{L}^1 = \mathcal{A}^1_{>0} / \mathcal{A}^1_{>0} \mathcal{A}^1_{>0} \).  Glanois showed the following.
	
	\begin{Lem}[{Motivic Galois descent, \cite[Corollary 5.1.3]{GlanoisThesis16}, \cite[Theorem 3.8]{GlanoisBasis16}}]\label{lem:desc}
		Let \( \mathfrak{z} \in \mathcal{H}^2 \) be a linear combination of motivic alternating MZV's.  Then \( \mathfrak{z} \) is a linear combination of motivic non-alternating MZV's if and only if \( D_1(\mathfrak{z}) = 0 \) and \( D_{2r+1}(\mathfrak{z}) \in \mathcal{L}^1 \otimes \mathcal{H}^1 \), for \( r \geq 1 \).
	\end{Lem}
	
	In \cite{ZhengProof}, the authors make progress towards a proof of \autoref{conj:xu} using Glanois' motivic descent criterion in \autoref{lem:desc}.  Unfortunately, this argument has a limitation.  Glanois' criterion gives no information about the \emph{depth} of the resulting non-alternating combination, nor does the motivic coaction allow one to fix the depth of classical multiple zeta values, as the coradical filtration and depth filtration do not agree \cite{BrownDepthGraded21}.  Without giving an explicit formula which could be checked motivically, \cite{ZhengProof} can -- at best -- only show \( \zeta(2,\overline{2k}) + 2k \cdot \zeta(1, \overline{2k+1}) \) is a linear combination of non-alternating MZV's. 
	
	\begin{Eg} 
		The following evaluation for \( \zeta(\overline{3}, \overline{9}) \) contains the (conjecturally) irreducible depth 4 non-alternating multiple zeta value \( \zeta(1,1,4,6) \).  This evaluation is verified by the Multiple Zeta Value Data Mine \cite{mzvDM} (hence is motivic); more precisely \cite[Eqn. (10.6)]{mzvDM} calls a version of this identity a ``pushdown'' of \( \zeta(1,1,4,6) \) to depth 2.  The identity reads
	\begin{equation}
	\label{eqn;pushdown}
	\begin{aligned}[c]
		\zeta(\overline{3},\overline{9})= & \frac{9}{64} \zeta (1,1,4,6)
	-\frac{371 \zeta (3,9)}{1024}
	-\frac{27}{64} \zeta (2) \zeta (3,7)
	-\frac{27}{128} \zeta (4) \zeta (3,5)
	\\
	& {}
	 +\frac{3131}{1024} \zeta (3) \zeta (9) 
	-\frac{321}{512} \zeta (5) \zeta (7)
	-\frac{3}{256} \zeta (3)^4 
	-\frac{45}{32} \zeta (2) \zeta (3) \zeta (7)
	\\
	&{}
	-\frac{63}{128} \zeta (2) \zeta (5)^2
	+\frac{9}{128} \zeta (4) \zeta (3) \zeta (5)
	+\frac{81}{256} \zeta (6) \zeta (3)^2
	+\frac{353139 \zeta (12)}{2830336} \,.
	\end{aligned}
	\end{equation}
	This gives a Galois descent of \( \zeta(\overline{3},\overline{9}) \) to level 1, but shows that in general such descents \emph{do not} preserve the depth.  (See similar discussion in Remark A.4 in \cite{CK2242}.)
	
	Certainly \( \zeta(\overline{3},\overline{9}) \) satisfies the refined descent criterion for double zeta values given in Lemma 2.1 \cite{ZhengProof}, in particular one can check \( D_1 \zeta^\mot(\overline{3},\overline{9}) = D_{11} \zeta^\mot(\overline{3},\overline{9}) = 0 \) directly (since we do have a Galois descent, we get this automatically in any case).  But the depth of the Galois descent necessarily increases, so we cannot possibly obtain a proof of Xu's Conjecture, using just this descent criterion.
	\end{Eg}
	\medskip
	
	In this note, we establish the following explicit identity proving Xu's Conjecture (\autoref{conj:xu}).
	\begin{Thm}\label{thm:1}
		For any integer \( k \geq 1 \), the following evaluation in terms of non-alternating MZV's of depth \( \leq 2 \) holds,
	\begin{align*}
		&\hspace{-1em} \zeta(2,\overline{2k}) + 2k \zeta(1,\overline{2k+1}) = 
	\\
	& 
	\begin{aligned}[t]
	& \zeta(2k,2)  - \sum_{i=2}^{2k} \frac{1}{2^i} \bigg\{ \!  \binom{i-1}{1} \zeta(2k+2-i,i) + \binom{i-1}{2k-1} \zeta(i, 2k+2-i) \! \bigg\} \\
	& {} - \sum_{r=2}^{2k} \bigg\{ \! (-1)^r (1 - 2^{-r}) \binom{r-1}{1} + \binom{r-1}{2k-1} (1 - 2^{1-r}) \! \bigg\} \zeta(r)\zeta(2k+2-r)  \\
	& {} + \big(2 + (2k-1)(1+2^{-2k-2})\big) \zeta(2k+2) \,.
	\end{aligned}
	\end{align*}
	This identity also holds on the motivic level.
	\end{Thm}

	\begin{Eg}
	When \( k = 5 \), we have the weight 12 evaluation, 
	\begin{align*}
	& \hspace{-1em} \zeta (2,\overline{10}) + 10 \zeta (1,\overline{11}) = \\
	& 
	-\frac{9}{1024} \zeta (2,10)
	-\frac{1}{64} \zeta (3,9)
	-\frac{7}{256} \zeta (4,8)
	-\frac{3}{64} \zeta(5,7)
	-\frac{5}{64} \zeta (6,6) 
	-\frac{1}{8} \zeta (7,5)
	\\
	&
	-\frac{3}{16} \zeta (8,4)
	-\frac{1}{4} \zeta (9,3)
	+\frac{767}{1024} \zeta(10,2)
	\,\, + \,\,\frac{623}{64} \zeta (3) \zeta (9)
	+\frac{629}{64} \zeta (5) \zeta (7)
	\\
	&
	-\frac{10997}{1024} \zeta (2) \zeta (10)
	-\frac{315}{64} \zeta (6)^2
	-\frac{2505}{256} \zeta (4) \zeta (8)
	+\frac{45065}{4096} \zeta (12) \,.
	\end{align*}
	This can still be verified with the multiple zeta value Data Mine \cite{mzvDM}, as the weight is sufficiently low.
	\end{Eg}

	\section{Proof of \protect\autoref{thm:1}}
	
	\Autoref{thm:1} follows by interpreting \( \zeta(2,\overline{2k}) + 2k \zeta(1,\overline{2k+1})\) as the regularisation of one particular \emph{regularised} alternating double zeta value \( \zeta_1^{\shuffle,0}(1, \overline{2k+1})\), and applying explicit identities on alternating double zeta values established in \cite{CK2242}.

	\subsection*{Iterated integrals and regularisation}  We briefly recall the setup of regularising iterated integral, and the shuffle regularisation of alternating MZV's.  The iterated integral is defined by
	\begin{equation}\label{eqn:int}
		I(a; x_1,\ldots,x_n, b) = \int_{a < t_1 < \cdots < t_n < b} \frac{\dd t_1}{t_1 - x_1} \cdots \wedge \frac{\dd t_n}{t_n - x_n} \,.
	\end{equation}
	The iterated integrals multiply with the shuffle product:
	\[
		I(a; x_1,\ldots,x_n, b) I(a; x_{n+1},\ldots,x_{n+m}, b) = \sum_{\sigma \in \Sigma_{n,m}} I(a; x_{\sigma(1)}, \ldots, x_{\sigma(n+m)}; b) \,,
	\]
	where \( \Sigma_{n,m} \) is the set of $(n,m)$-shuffles, i.e. permutations \( \sigma \in \mathfrak{S}_{n+m} \), such that \( \sigma^{-1}(1) < \cdots < \sigma^{-1}(n) \) and \( \sigma^{-1}(n+1) < \cdots < \sigma^{-1}(n+m) \).
	When \( x_1 =  a \) or \( x_n = b \), the integral diverges; by considering the asymptotic expansion of \( I(a+\epsilon; x_1, \ldots, x_n; b-\epsilon) \) as a polynomial in \( \log(\epsilon) \), one defines the regularised value as the constant term (cf. `canonical regularisation' \cite[\S9]{GoncharovMultiple01}); write \( I^{\shuffle} \) for the regularised value.  The regularised integrals also satisfy the shuffle product.  The regularisation of a convergent integral is simply the original value of that integral.

	For example,
	\begin{equation}\label{eqn:reg0}
		I^{\shuffle}(0; 0; 1) = [\log(\epsilon)^0] \int_{0+\epsilon}^{1-\epsilon} \frac{\dd t}{t} 
		 = [\log(\epsilon)^0] ( \log(1 - \epsilon) - \log(0 + \epsilon) ) 
		 = \log(1) = 0
	\end{equation}
	
	We recall that alternating MZV's have an iterated integral representation
	\[
		\zeta(k_1,\ldots,k_d; \epsilon_1,\ldots,\epsilon_d) = (-1)^d I(0; {\overbrace{\eta_1, 0, \ldots, 0}^{k_1}}, \ldots, {\overbrace{\eta_d, 0, \ldots, 0}^{k_d}};  1) \,, \qquad \eta_j = \textstyle\prod_{i=j}^d \epsilon_i \,.
	\]
	We define shuffle regularised MZV's as
	\[
		\zeta_{k_0}^{\shuffle}(k_1,\ldots,k_d; \epsilon_1,\ldots,\epsilon_d) = (-1)^d I(0; \smash{\overbrace{0,\ldots,0}^{k_0}}, \smash{\overbrace{\eta_1, 0, \ldots, 0}^{k_1}}, \ldots, \smash{\overbrace{\eta_d, 0, \ldots, 0}^{k_d}};  1) \,, \qquad \eta_j = \textstyle\prod_{i=j}^d \epsilon_i \,,
	\]
	allowing an arbitrary number of 0's at the start of the integral.  (This approach also defines a regularisation for \( \zeta(k_1,\ldots,k_d; \epsilon_1,\ldots,\epsilon) \), with \( (k_d,\epsilon_d) = (1,1) \).)  From \eqref{eqn:reg0} (and a similar computation for \( \zeta^\shuffle(1) = 0 \), we have obtained the shuffle regularisation of alternating MZV's, with regularisation \( \zeta_1^{\shuffle}(\emptyset) = \zeta^{\shuffle}(1) = 0 \).

	\subsection*{Regularisation expression}
	Write \( \{0\}^n = 0,\ldots,0 \), with \( n \) repetitions of \( 0 \).  By the shuffle product of regularised integrals, and \eqref{eqn:reg0}, we have
	\begin{equation}\label{eqn:goal} 
	\begin{aligned}[c]
		 \zeta(2,\overline{2k}) + 2k \cdot \zeta(1, \overline{2k+1}) 
		& = I(0; -1, 0, -1, \{0\}^{2k-1}; 1) + 2k \cdot I(0; -1, -1, \{0\}^{2k}; 1) \\
		& = I^{\shuffle}(0; 0, 1) I(0; -1, -1, \{0\}^{2k-1}; 1) - I^{\shuffle}(0; 0, -1, -1, \{0\}^{2k-1}; 1) \\
		& = -\zeta_1^{\shuffle,0}(1, \overline{2k+1})
	\end{aligned}
	\end{equation}
	
	\subsection*{Dihedral symmetry}
	We recall the following identity -- an explicit version of Glanois' dihedral symmetry \cite[Corollary 4.2.6]{GlanoisThesis16} -- given in Equation (A.8) \cite{CK2242} (we have interchanged \( k \leftrightarrow \ell \) for convenience),
	\begin{equation}\label{eqn:dihedral}
	\begin{aligned}
		& \zeta^{\shuffle}_{2\ell-1}(1,\overline{2k}) - \zeta(\overline{2k},\overline{2\ell}) 
		\\
		& =
		\binom{2k+2\ell-1}{2\ell-1} \zeta(\overline{2k+2\ell}) - \sum_{r=1}^{2k+2\ell-2} \bigg( \! (-1)^r \binom{r-1}{2\ell-1} + \binom{r-1}{2k-1} \! \bigg) \zeta(\overline{r}) \zeta(2k+2\ell-r) \,.
	\end{aligned}
	\end{equation}
	Notice that the \( r = 1 \) term of the sum vanishes, so \( \zeta(\overline{1}) \) does not enter the result.  
		
	This identity was established as follows: Regularising \( \zeta^{\shuffle}_{z-1}(\overline{\alpha}, \beta) + \zeta^{\shuffle}_{z-1}(\beta, \overline{\alpha}) \) and applying the stuffle antipode  reduces the result to products of depth 1 alternating MZV's.  Likewise, the shuffle antipode expresses \( \zeta^{\shuffle}_{z-1}(\overline{\alpha}, \beta) + (-1)^{z-1 + \alpha+\beta} \zeta^{\shuffle}_{\beta-1}(\overline{\alpha}, \overline{z}) \) as products of depth 1 alternating MZV's.  Taking the difference in the case \( (\alpha,\beta,z) = (2k,1,2\ell) \) gives the dihedral symmetry in \eqref{eqn:dihedral}.
	
	Specialising \eqref{eqn:dihedral} to \( \ell = 1 \) gives
	\begin{equation}\label{eqn:dihedral2k2}
	\begin{aligned}[c]
		& \zeta^{\shuffle}_{1}(1,\overline{2k}) - \zeta(\overline{2k},\overline{2}) 
		\\
		& =
		\binom{2k+1}{1} \zeta(\overline{2k+2}) - \sum_{r=2}^{2k} \bigg( \! (-1)^r \binom{r-1}{1} + \binom{r-1}{2k-1} \! \bigg) \zeta(\overline{r}) \zeta(2k+2-r) \,.
	\end{aligned}
	\end{equation}
	\medskip
		
	\subsection*{Descent} Next, we recall the explicit Galois descent of \( \zeta(\overline{2k},\overline{2\ell}) \) to non-alternating double zeta values established in Proposition A.3 \cite{CK2242} (again we have interchanged \( k \leftrightarrow \ell \) for convenience),
\begin{equation}
\label{eqn:zaltevev}
\begin{aligned}[c]
\zeta(\overline{2k},\overline{2\ell}) = {}
& \sum_{i=2}^{2k+2\ell - 2} 2^{-i} \bigg\{ \! \binom{i-1}{2\ell-1} \zeta(2k+2\ell-i,i) + \binom{i-1}{2k-1} \zeta(i,2k+2\ell-i) \! \bigg\} \\
& -\zeta(2k,2\ell) + \sum_{r=2}^{2k+2\ell-2} (-2)^{-r} \binom{r-1}{2\ell-1}  \zeta(r) \zeta(2k+2\ell-r) \\
& - 2^{-2k-2\ell} \bigg\{ \! 2 \binom{2k+2\ell-2}{2\ell-1} + \binom{2k+2\ell-1}{2\ell-1} \! \bigg\} \zeta(2k+2\ell) \,.
\end{aligned}
\end{equation}

	This identity was established by solving a system of simultaneous equations coming from the dihedral symmetry \eqref{eqn:dihedral} and a generalised doubling relation \cite[\S4]{mzvDM} (see \cite[\S 14.2.5]{ZhaoBook} also).  The former  expressed \( \zeta^{\shuffle}_{2\ell-1}(1,\overline{2k}) - \zeta(\overline{2k},\overline{2\ell}) \) via products of depth 1 MZV's, while the latter expressed \( \zeta^{\shuffle}_{2\ell-1}(1,\overline{2k}) + \zeta(\overline{2k},\overline{2\ell}) \) via non-alternating double zeta values and products of depth 1 MZV's.
	
	Setting \( \ell = 2 \) in \eqref{eqn:zaltevev} produces the following expression for \( \zeta(\overline{2k},\overline{2}) \),
\begin{equation}
\label{eqn:zalt2k2}
\begin{aligned}[c]
\zeta(\overline{2k},\overline{2}) = {}
 & \sum_{i=2}^{2k} \frac{1}{2^i} \bigg\{ \! \binom{i-1}{1} \zeta(2k+2-i,i) + \binom{i-1}{2k-1} \zeta(i,2k+2-i) \! \bigg\} \\
& -\zeta(2k,2) + \sum_{r=2}^{2k} \frac{r-1}{(-2)^r}  \zeta(r) \zeta(2k+2-r) - \frac{6k+1}{2^{2k+2}} \zeta(2k+2) \,.
\end{aligned}
\end{equation}

	\subsection*{Combining the results}
	It is now simple to substitute \eqref{eqn:zalt2k2} into \eqref{eqn:dihedral2k2} to obtain an expression for \( \zeta_1^{\shuffle,0}(1, \overline{2k}) \) in terms of depth 2 non-alternating MZV's.  Using this evaluation for \(  \zeta_1^{\shuffle,0}(1, \overline{2k})  \) in \eqref{eqn:goal} gives an evaluation for \( \zeta(2,\overline{2k}) + 2k \zeta(1, \overline{2k+1})  \).  This proves Xu's Conjecture (\autoref{conj:xu}).
	
	The explicit expression in \autoref{thm:1} follows by some straight forward identities and simplifications.  In particular we have: 
	\begin{enumerate}
		\item\label{enum1}  written depth 1 alternating MZV's as non-alternating MZV's, via \( \zeta(\overline{n}) = -(1-2^{1-n})\zeta(n) \), for \( n \geq 2 \) (as noted, no \( \zeta(\overline{1}) \) terms are present); 
		\item combined and simplified the \( \zeta(2k+2) \) terms, as well as the sums over products of single zeta values in \eqref{eqn:dihedral2k2} and \eqref{eqn:zalt2k2}.
	\end{enumerate}
	This complete the proof of \autoref{thm:1}.  Recall that both \eqref{eqn:dihedral} and \eqref{eqn:zaltevev} were verified to be motivic in \cite[Appendix B]{CK2242}, hence \autoref{thm:1} also holds on the motivic level. \hfill \qedsymbol

	\bibliographystyle{habbrv2} 
	\bibliography{bibliography.bib}

\end{document}